\documentclass[12pt]{amsart}
\usepackage{a4wide}
\usepackage[latin1, utf8]{inputenc}

\usepackage[hang,flushmargin]{footmisc}

\usepackage{amsmath,amssymb,amsfonts}
\usepackage[foot]{amsaddr}
\usepackage{graphicx}
\usepackage{hyperref}
\usepackage{orcidlink}
\usepackage{doi}

\usepackage[font=footnotesize,labelfont=footnotesize]{caption}

\usepackage[square,numbers]{natbib}

\usepackage{mathtools}
\mathtoolsset{showonlyrefs}

\numberwithin{equation}{section}
\newtheorem{theorem}{Theorem}[section]

\newcommand{\R}{{\mathbb R}}

\newcommand{\Z}{{\mathbb Z}}

\newcommand{\vol}{{\rm vol}}

\usepackage{wasysym}

\author[M.~Faulhuber]{Markus Faulhuber \orcidlink{0000-0002-7576-5724}}
\address{Faculty of Mathematics, University of Vienna \newline Oskar-Morgenstern-Platz 1, 1090 Vienna, Austria}
\email{markus.faulhuber@univie.ac.at}

\begin{document}
\title[Polarization for the honeycomb]{The polarization problem for the honeycomb structure}

\thanks{This research was funded in whole or in part by the Austrian Science Fund (FWF) [\href{https://doi.org/10.55776/P33217}{10.55776/PAT5102224}].}



\begin{abstract}
    We study the polarization problem in dimension 2 for the honeycomb structure and compare it to the maximal polarization lattice, the hexagonal lattice. As expected, the hexagonal lattice has higher polarization than the honeycomb at all densities.
\end{abstract}

\vspace*{-\baselineskip}
\maketitle

\section{Introduction and notation}
The \textit{polarization problem} has its origin in potential analysis and the study of Chebyshev constants and the work of Fekete \cite{Fek23} and Poly\'a and Szeg\H{o} \cite{PolSze31}. It has been intensively studied in recent years for the Riesz potential on the ($d$-dimensional) sphere, see, e.g., \cite{AmbBalErd12}, \cite{Bor22}, \cite{BorHarRezSaf18}, \cite{BoyDraHarSafSto23}, \cite{HarKenSaf13}, to name just a few articles. We consider the problem for periodic point configurations in the plane and for exponential potentials of squared distance, i.e., Gaussians. A problem for potentials $f(r^2)$ with $f:\R_+ \to \R_+$, being decreasing and convex was formulated for lattices by Hardin, Petrache and Saff \cite{HarPetSaf22} and answered by B\'etermin, Faulhuber, and Steinerberger \cite{BetFauSte21}. Most recently, Bachoc, Moustrou, Vallentin, and Zimmermann characterized \textit{cold spots} of many root lattices \cite{BacMouValZim25}, which is a first step towards tackling the polarization problem in $\R^d$.

We say a function $f: \R_+ \to \R$, $r \mapsto f(r)$ is completely monotone if $(-1)^k f^{(k)}(r) \geq 0$. We will consider the following problem:
\begin{equation}\label{eq:polarization}
    \max_{\Gamma \subset \R^d} \min_{z \in \R^d} \sum_{\gamma \in \Gamma} f(|\gamma-z|^2).
\end{equation}
This is known as the polarization problem in Euclidean space $\R^d$. The minimum depends generally on the geometry and density (defined below) of $\Gamma$, and the potential $f$. The maximizer among (periodic) point configurations $\Gamma$ is sought for a fixed density of $\Gamma$. In rare cases, a maximizer may be optimal among all densities (scales), in which case it becomes a universal polarization maximizer, much like universally optimal energy minimizers in the sense of Cohn and Kumar \cite{CohKum07}. For energy minimization, solutions have been derived in dimensions 8 and 24 by Cohn, Kumar, Miller, Radchenko, and Viazovska \cite{Coh-Via22}. The energy minimization problem in dimension 2 was solved among lattices by Montgomery~\cite{Mon88}. For polarization, we only know that the hexagonal lattice $\Lambda_2$ (the $\mathsf{A}_2$ root lattice) is a universal polarization maximizer among lattices (of any fixed density) and completely monotone potentials of squared distance \cite{BetFauSte21}.
In this note we prove the following polarization result.

\begin{theorem}\label{thm:main}
    Denote the hexagonal lattice by $\Lambda_2$ and the honeycomb structure by $\Gamma_2$, both of density 1. For an arbitrary but fixed density $\rho > 0$ and any completely monotone potential $f: \R_+ \to \R$ (of sufficient decay at infinity) of squared distance, the hexagonal lattice of density $\rho$ has higher polarization than the honeycomb structure of density $\rho$:
    \begin{equation}
        \min_{z} \sum_{\lambda \in \rho^{-1/2}\Lambda_2} f(|\lambda-z|^2) > \min_{z} \sum_{\lambda \in \rho^{-1/2}\Gamma_2} f(|\gamma-z|^2).
    \end{equation}
\end{theorem}

By the Bernstein-Widder theorem \cite{Ber28}, \cite{Wid41}, a completely monotone function of squared distance can be written as
\begin{equation}\label{eq:Bernstein-Widder}
    f(r^2) = \int_{\R_+} e^{-\alpha r^2} \, d \mu_f(\alpha),
\end{equation}
for a non-negative Borel measure $\mu_f$. In accordance with \cite{BacMouValZim25}, we say a structure $\Gamma$ has a cold spot in $z_0$ if (fix $f$ for the moment) $z_0$ is a minimizer in \eqref{eq:polarization}. It is a universal cold spot if it does not depends on the density of $\Gamma$. For the potential $\phi_\alpha(r^2) = e^{-\pi \alpha r^2}$, if a $\Gamma$ solve \eqref{eq:polarization} for any fixed density (which corresponds to being a maximizer at density 1 and all scales $\alpha > 0$ of the Gaussian), then it is a universal polarization maximizer for all completely monotone functions of squared distance $f(r^2)$ (given fast enough decay of $f$). In this case the universal cold spot is also independent of $f$. Yet, the existence of a universally polarization maximizers is open, except for lattices in dimension 2 \cite{BetFauSte21}.

A (full-rank) lattice $\Lambda$ in $\R^d$ is a discrete co-compact subgroup (see \cite{ConSlo98} for a full introduction) and its density is the average number of lattice points per volume, given by
\begin{equation}
    \rho(\Lambda) = 1/\vol(\R^d/\Lambda).
\end{equation}
A periodic configuration is the union of relatively shifted copies of a lattice $\Lambda \subset \R^d$;
\begin{equation}
    \Gamma = \bigcup_{n=1}^N (\Lambda + x_n), \quad x_m-x_n \notin \Lambda, m \neq n.
\end{equation}
The density of $\Gamma$ is then given by $\rho(\Gamma) = N/\vol(\R^d/\Lambda)$ (cf.\ \cite[Sec.~9]{CohKum07}).

\section{Universal cold spots and the honeycomb structure}
\begin{figure}[ht]
    \centering
    \includegraphics[width=0.75\linewidth]{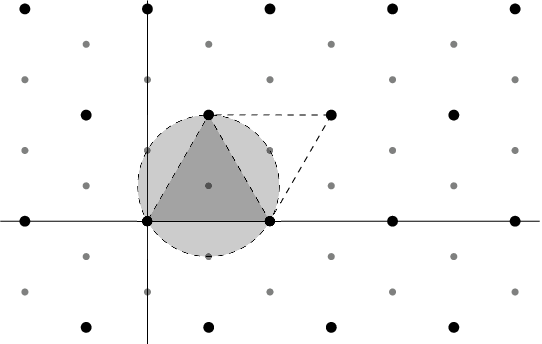}
    \captionsetup{width=\linewidth}
    \caption{A part of the hexagonal lattice $\Lambda_2$ with its fundamental cell, a fundamental triangle, and a covering circle centered at a deep hole. The large black dots are the lattice points and the smaller gray dots are deep holes.}
    \label{fig:A2_lattice}
\end{figure}

The honeycomb $\Gamma_2$ structure is the union of two hexagonal lattices, relatively shifted by a deep hole $z_0$. A deep hole of a lattice is a point in $\R^d$ maximizing the distance to the closest lattice point(s). For the hexagonal lattice, it is the circumcenter of the equilateral (fundamental) triangle formed by 3 lattice points (see Fig.~\ref{fig:A2_lattice}).

The hexagonal lattice of density 1 is, up to rotational equivalence, explicitly given by
\begin{equation}
    \Lambda_2 = \sqrt{\frac{1}{\sin(\pi/3)}}
    \begin{pmatrix}
        1 & \cos(\pi/3)\\
        0 & \sin(\pi/3)
    \end{pmatrix} \Z^2
    = \sqrt{\tfrac{2}{\sqrt{3}}}
    \begin{pmatrix}
        1 & \frac{1}{2}\\
        0 & \frac{\sqrt{3}}{2}
    \end{pmatrix} \Z^2
\end{equation}
The honeycomb structure of unity density is given by
\begin{equation}\label{eq:honeycomb}
    \Gamma_2 = \sqrt{2} \left( \Lambda_2 \cup (\Lambda_2+z_0) \right).
\end{equation}
The factor $\sqrt{2}$ scales $\Gamma_2$ appropriately and we always write $z_0$ for the deep hole(s) of $\Lambda_2$.

For the hexagonal lattice $\Lambda_2$ and the potential $\phi_\alpha(r) = e^{-\pi \alpha r}$ (the $\pi$ in the exponent is for technical reasons, due to our normalization \eqref{eq:FT} of the Fourier transform) we have
\begin{equation}\label{eq:Lambda2_z0}
    \sum_{\lambda \in \Lambda_2} \phi_\alpha(|\lambda-z_0|^2) \leq \sum_{\lambda \in \Lambda_2} \phi_\alpha(|\lambda-z|^2), \quad \forall \alpha > 0, \ z \in \R^2.
\end{equation}
This result is due to Baernstein \cite{Baernstein_HeatKernel_1997} and shows that deep holes are universal cold spots of the hexagonal lattice. Note that if $z_0 \in \R^2/\Lambda_2$ is a deep hole, then $2 z_0 \notin \Lambda_2+z_0$ is also a deep hole. These two deep holes are inequivalent; all other deep holes can be reached from these by translation by $\Lambda_2$. For the honeycomb structure $\Gamma_2$, we have, by setting $\widetilde{z} = \sqrt{2} z$,
\begin{align}
    \sum_{\gamma \in \Gamma_2} \phi_\alpha(|\gamma - \widetilde{z}|^2)
    & = \sum_{\lambda \in \Lambda_2} \phi_\alpha(2|\lambda-z|^2) + \sum_{\lambda \in \Lambda_2} \phi_\alpha(2|\lambda+z_0-z|^2)\\
    & = \sum_{\lambda \in \Lambda_2} \phi_{2\alpha}(|\lambda-z|^2) + \sum_{\lambda \in \Lambda_2} \phi_{2\alpha}(|\lambda+z_0-z|^2).\label{eq:honeycomb_split}
\end{align}
Note that scaling the arguments is necessary, due to \eqref{eq:honeycomb}. The scaling is put into the parameter $\alpha$ of the Gaussian $\phi_\alpha(r^2)$. This will also happen in sequential computations.

As $2 z_0$ is a deep hole of $\Lambda_2$ and of $\Lambda_2+z_0$ (corresponding to $z_0$ of $\Lambda_2$) we have by \eqref{eq:Lambda2_z0}
\begin{equation}
    \sum_{\gamma \in \Gamma_2} \phi_\alpha(|\gamma - 2 \widetilde{z}_0|^2) \leq \sum_{\gamma \in \Gamma_2} \phi_\alpha(|\gamma - z|^2), \quad \forall \alpha > 0, \ z \in \R^2,
\end{equation}
where $\widetilde{z}_0 = \sqrt{2} z_0$ is a deep hole of the scaled hexagonal lattice $\sqrt{2} \Lambda_2$.
\begin{figure}
    \centering
    \includegraphics[width=0.75\linewidth]{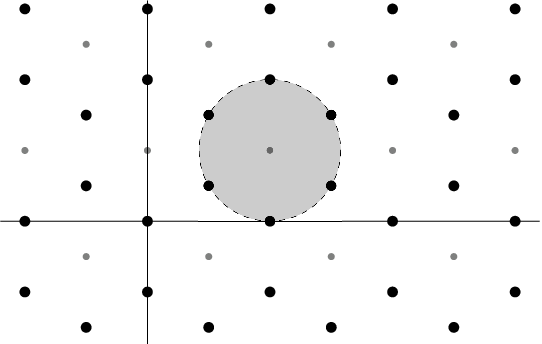}
    \captionsetup{width=\linewidth}
    \caption{The honeycomb structure $\Gamma_2$ is a union of the hexagonal lattice and a copy shifted by a deep hole. The second, inequivalent deep hole of $\Lambda_2$ is a deep hole of $\Gamma_2$. The big black dots are points of $\Gamma_2$ while the smaller gray dots are its deep holes.}
    \label{fig:honeycomb}
\end{figure}

\section{Proof or Theorem \ref{thm:main}}
Most methods we use are basic, the Poisson summation formula being already an advanced technique. However, we build on deep insights by Baernstein \cite{Baernstein_HeatKernel_1997} on the universality of cold spots in the hexagonal lattice. At one point we will use a result of Ramanujan (see \cite[Chap.~33]{RamanujanV}), put on solid mathematical ground by Borwein and Borwein~\cite{BorBor_Cubic_91}.

By \eqref{eq:Bernstein-Widder} it suffices to consider the potential $\phi_\alpha(|r^2|)$ and to show that
\begin{equation}\label{eq:goal}
    \sum_{\gamma \in \Gamma_2} \phi_{\alpha}(|\gamma - 2 \widetilde{z}_0|^2) <  \sum_{\lambda \in \Lambda_2} \phi_\alpha(|\lambda-z_0|^2), \quad \forall \alpha > 0,
\end{equation}
in order to prove Theorem~\ref{thm:main}. We start with noting that
\begin{equation}
    \sum_{\lambda \in \Lambda_2} \phi_\alpha(|\lambda-2z_0|^2) = \sum_{\lambda \in \Lambda_2} \phi_\alpha(|\lambda-z_0|^2),
\end{equation}
as both $z_0$ and $2z_0$ are deep holes of $\Lambda_2$ and, thus, universal cold spots \cite{Baernstein_HeatKernel_1997}. We rewrite the entities from \eqref{eq:goal} over the integer lattice $\Z^2$ by a linear coordinate change (cf.\ \cite{BetFauSte21}, \cite{Mon88}):
\begin{align}
    \sum_{\lambda \in \Lambda_2} \phi_\alpha(|\lambda-z_0|^2)
    & = \sum_{(k,l) \in \Z^2} e^{- \frac{2 \pi \alpha}{\sqrt{3}} ((k-1/3)^2 + (k-1/3)(l-1/3)+(l-1/3)^2)},
    \\
    \sum_{\gamma \in \Gamma_2} \phi_{\alpha}(|\gamma - 2 \widetilde{z}_0|^2)
    & = \sum_{\lambda \in \Lambda_2} \phi_{2\alpha}(|\lambda-2 z_0|^2) + \sum_{\lambda \in \Lambda_2} \phi_{2\alpha}(|\lambda-z_0|^2) = 2 \sum_{\lambda \in \Lambda_2} \phi_{2\alpha}(|\lambda-z_0|^2)\\
    & = 2 \sum_{(k,l) \in \Z^2} e^{- \frac{4 \pi \alpha}{\sqrt{3}} ((k-1/3)^2 + (k-1/3)(l-1/3)+(l-1/3)^2)}.
\end{align}
To split the series for the honeycomb we used \eqref{eq:honeycomb_split}. It is not hard to observe that
\begin{equation}
    \min_{(k,l) \in \Z^2} \left((k-1/3)^2 + (k-1/3)(l-1/3)+(l-1/3)^2\right) = 1/3,
\end{equation}
which is achieved if and only if
\begin{equation}
    (k,l) \in \{(0,0),(1,0),(0,1)\}.
\end{equation}
As a next step, we note that for $x \in \R_+$
\begin{equation}
    e^{-x} > 2 e^{-2x}
    \quad \Longleftrightarrow \quad
    x > \log(2).
\end{equation}
This particularly implies that for all $(k,l) \in \Z^2$
\begin{align}
    & e^{- \frac{2 \pi \alpha}{\sqrt{3}} ((k-1/3)^2 + (k-1/3)(l-1/3)+(l-1/3)^2)} > e^{- \frac{4 \pi \alpha}{\sqrt{3}} ((k-1/3)^2 + (k-1/3)(l-1/3)+(l-1/3)^2)}\\
    \quad \Longleftrightarrow \quad
    & \frac{2 \pi \alpha}{3 \sqrt{3}} > \log(2)\\
    \quad \Longleftrightarrow \quad
    & \alpha > \frac{3 \sqrt{3} \log(2)}{2 \pi} \approx 0.573228.
\end{align}
It readily follows that \eqref{eq:goal} holds for $\alpha$ large enough:
\begin{equation}
    \sum_{\lambda \in \Lambda_2} \phi_\alpha(|\lambda-z_0|^2) > \sum_{\gamma \in \Gamma_2} \phi_\alpha(|\gamma-\widetilde{z}_0|^2), \quad \forall \alpha > \frac{3 \sqrt{3} \log(2)}{2 \pi}.
\end{equation}
To show that the polarization of $\Lambda_2$ is larger than that of $\Gamma_2$ on small scales $\alpha$ (equivalently for high densities), we employ the Poisson summation formula as usual in such situations.

We recall that the dual lattice $\Lambda^*$ is characterized by $\Lambda^* = \{ \lambda^* \in \R^d \mid \lambda \cdot \lambda^* \in \Z \ \forall \lambda \in \Lambda \}$, the dot $\cdot$ denoting the Euclidean inner product. We denote the Fourier transform of $f$ by
\begin{equation}\label{eq:FT}
    \widehat{f}(y) = \int_{\R^d} f(x) e^{-2 \pi i x \cdot y} \, dx.
\end{equation}
The Poisson summation formula, in this case, reads
\begin{equation}
    \sum_{\lambda \in \Lambda} f(\lambda+z) = \frac{1}{\vol(\R^d/\Lambda)} \sum_{\lambda^* \in \Lambda^*} \widehat{f}(\lambda^*) e^{2 \pi i \lambda^* \cdot z}, \quad z \in \R^d.
\end{equation}

Our situation is special; any lattice $\Lambda \subset \R^2$ (of density 1) is self-dual. More precisely, if $\Lambda \subset \R^2$, $\vol(\R^2/\Lambda) = 1$, then $\Lambda^* = J \Lambda$ (cf.\ \cite[App.~A]{BetFauSte21}),
where
\begin{equation}
    J = \begin{pmatrix}
        0 & 1\\
        -1 & 0
    \end{pmatrix}
\end{equation}
is a rotation by 90 degrees. Generally, the Gaussian $x \mapsto e^{-\pi |x|^2}$, $x \in \R^d$, is an eigenfunction of the Fourier transform. In particular, we have \cite[App.~A]{Fol89}
\begin{equation}
    \widehat{\phi_\alpha}(y^2) = \alpha^{-d/2} \phi_{1/\alpha}(y^2), \quad y \in \R^d.
\end{equation}
\enlargethispage{\baselineskip}
So, in our very particular situation, we have (note that $-\lambda \in \Lambda$ for any lattice)
\begin{equation}
    \sum_{\lambda \in \Lambda_2} \phi_{\alpha}(|\lambda-z_0|^2) = \sum_{\lambda \in \Lambda_2} \phi_{\alpha}(|\lambda+z_0|^2) = \alpha^{-1} \sum_{\lambda \in \Lambda_2} \phi_{1/\alpha}(|\lambda|^2) e^{2 \pi i (J\lambda) \cdot z_0}.
\end{equation}
We used the symmetry of $\Lambda_2$, the Poisson summation formula, the fact that $\vol(\R^2/\Lambda_2)=1$, that $\Lambda_2^* = J \Lambda_2$, and that $\phi_\alpha(|r|^2)$ is rotational invariant. Note that the left side is real-valued, so the right side must also be real-valued. With the same arguments, we obtain a dual formula for the polarization of $\Gamma_2$:
\begin{equation}
     \sum_{\gamma \in \Gamma_2} \phi_{\alpha}(|\gamma - 2 \widetilde{z}_0|^2)
    = 2 \sum_{\lambda \in \Lambda_2} \phi_{2\alpha}(|\lambda+z_0|^2) = \alpha^{-1} \sum_{\lambda \in \Lambda_2} \phi_{1/(2\alpha)}(|\lambda|^2) e^{2 \pi i (J\lambda) \cdot z_0}.
\end{equation}
Thus, by making the substitution $\alpha \mapsto \alpha^{-1}$, we seek to show that
\begin{equation}
    \sum_{\lambda \in \Lambda_2} \phi_{\alpha}(|\lambda|^2) e^{2 \pi i (J\lambda) \cdot z_0} > \sum_{\lambda \in \Lambda_2} \phi_{\alpha/2}(|\lambda|^2) e^{2 \pi i (J\lambda) \cdot z_0}, \quad \alpha >  \frac{2 \pi}{3 \sqrt{3} \log(2)} \approx 1.74451.
\end{equation}
By splitting the complex exponential into its real and imaginary parts and using the linear coordinate change from before, we actually obtain (cf.\ \cite{BetKnu18-Born})
\begin{align}
    \sum_{\lambda \in \Lambda_2} \phi_{\alpha}(|\lambda|^2) e^{2 \pi i (J\lambda) \cdot z_0}
    & = \sum_{\lambda \in \Lambda_2} \phi_{\alpha}(|\lambda|^2) \cos(2 \pi (J\lambda) \cdot z_0)\\
    & = \sum_{(k,l)\in \Z^2} e^{- \frac{2 \pi \alpha}{\sqrt{3}} (k^2+kl+l^2)} \cos\left(2 \pi (k-l)/3 \right)
\end{align}
Note that $\cos(2\pi (k-l)/3) \in \{-1/2, 1\}$ for all $(k,l) \in \Z^2$. Also, the way in which the values are taken is determined by the quadratic form $q(k,l) = k^2+kl+l^2$ and
\begin{equation}
    \lim_{R \to \infty}\frac{|\{(k,l) \in \Z^2 \mid \cos(2\pi (k-l)/3)  = -1/2\} \cap [-R,R]^2|}{|\{(k,l) \in \Z^2 \mid \cos(2\pi (k-l)/3)  = 1\} \cap [-R,R]^2|} = 2.
\end{equation}

\begin{figure}[ht]
    \centering
    \includegraphics[width=0.465\linewidth]{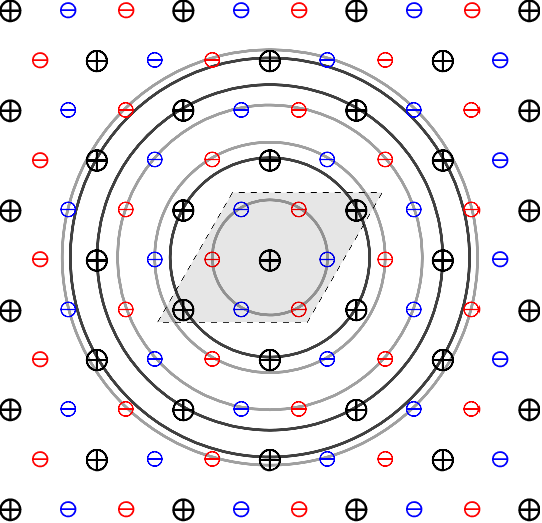}
    \hfill
    \includegraphics[width=0.45\linewidth]{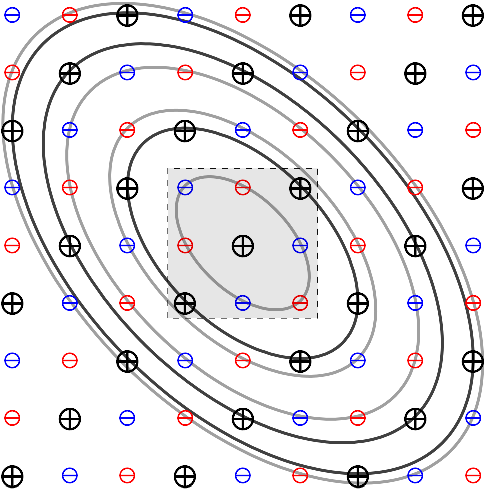}
    \captionsetup{width=\linewidth}
    \caption{Left: After performing a Poisson summation the polarization problem becomes a question on charged ionic crystals. The hexagonal lattice can then be decomposed into 3 (shifted) hexagonal lattices, each of third the density, 2 having negative charges of half weights and one having positive chagres.\\ Right: After a linear change of coordinates, the hexagonal setup is transferred to an anisotropic interaction problem on $\Z^2$. The anisotropic interaction is determined via the quadratic form $q(k,l)=k^2+kl+l^2$. In the square $[-1,1]^2$ we find twice as many negative charges of half weights as positive ones.}
    \label{fig:charges}
\end{figure}

So we have twice as many negative signs with weight 1/2 in the series as positive signs (with weight 1), leading to an overall ``neutral charge'' of $\Lambda_2$ (cf.\ \cite{BetFauKnu21}). In fact, we can decompose $\Lambda_2$ into 3 relatively shifted copies (shifted by deep holes), 2 copies having negative signs with half weights, and 1 copy with positive signs (and full weights) (cf.\ \cite{BetKnu18-Born}). This is the geometric interpretation (cf.\ \cite{FauGumSha25}) of the \textit{arithmetic part} of the cubic AGM of Browein and Borwein \cite{BorBor_Cubic_91}: for a complex number $q$ with $|q|<1$ set
\begin{equation}
    a(q) = \sum_{(k,l) \in \Z^2} q^{k^2+kl+l^2}
    \quad \text{ and } \quad
    b(q) = \sum_{(k,l) \in \Z^2} q^{k^2+kl+l^2} e^{2 \pi i \frac{k-l}{3}}.
\end{equation}
Then
\begin{equation}
    3 a(q^3) = a(q) + 2b(q).
\end{equation}
Setting $q = e^{-\frac{2 \pi \alpha}{\sqrt{3}}} < 1$, $\alpha > 0$, we have
\begin{equation}
    \sum_{\lambda \in \Lambda_2} \phi_{\alpha}(|\lambda|^2) e^{2 \pi i (J\lambda) \cdot z_0} = b(q) = \frac{3 a(q^3)-a(q)}{2},
\end{equation}
and
\begin{equation}
    \sum_{\lambda \in \Lambda_2} \phi_{\alpha/2}(|\lambda|^2) e^{2 \pi (J\lambda) \cdot z_0} = b(q^{1/2}) = \frac{3a(q^{3/2})-a(q^{1/2})}{2}.
\end{equation}
Similar to before, we make an auxiliary computation. We seek to find $x \in \R_+$ such that
\begin{equation}
    3 e^{-3x}- e^{-x} > 3e^{-\frac{3x}{2}} - e^{-\frac{x}{2}}.
\end{equation}
Substituting $s = e^{-x/2}$, this leads to finding $s \in [0,1]$ such that
\begin{equation}
    3s^6-3s^3-s^2+s = s \ (3s^5-3s^2-s+1) = s \  g(s) >0.
\end{equation}
Clearly, the value of $g(s)$ at $s=0$ is $g(0)=1$. Moreover, $g(1/5) = 2128/15625 > 0$. We check the sign of the derivative. We have
\begin{equation}
    g'(s) = 15 s^4 - 6s - 1.
\end{equation}
For $0 < s < 1/5$ we have $15s^4 = 3/125 < 1$ and so $g'(s) < 0$ on $(0,1/5)$. Thus $g(s)$ is positive on the interval $(0,1/5)$ and, ultimately,
\begin{equation}
    s \ g(s) > 0, \ s \in (0,1/5)
    \quad \Longleftrightarrow \quad
    3 e^{-3x}- e^{-x} > 3e^{-\frac{3x}{2}} - e^{-\frac{x}{2}}, \ x \in (2 \log(5), \infty).
\end{equation}
We note that our estimates are far from optimal, but they imply
\begin{equation}
    \frac{3a(q^3)-a(q)}{2} > \frac{3a(q^{3/2})-a(q^{1/2})}{2}, \quad q = e^{-\frac{2 \pi \alpha}{\sqrt{3}}}, \quad \alpha > \frac{2 \log(5) \sqrt{3}}{2 \pi} \approx 0.88733,
\end{equation}
by comparing for all $(k,l) \in \Z^2 \setminus\{(0,0)\}$. Note that the value at the origin adds 1 to both sides of the inequality. Thus, the proof of Theorem~\ref{thm:main} is finished.
\hfill $\square$


\begin{thebibliography}{23}
\providecommand{\natexlab}[1]{#1}
\providecommand{\url}[1]{\texttt{#1}}
\expandafter\ifx\csname urlstyle\endcsname\relax
  \providecommand{\doi}[1]{doi: #1}\else
  \providecommand{\doi}{doi: \begingroup \urlstyle{rm}\Url}\fi

\bibitem{AmbBalErd12}
G.~Ambrus, K.~M. Ball, and T.~Erd\'elyi.
\newblock {Chebyshev constants for the unit circle}.
\newblock \emph{Bulletin of the London Mathematical Society}, 45:\penalty0
  236--248, 2012.
\newblock \doi{10.1112/blms/bds082}.

\bibitem{BacMouValZim25}
C.~Bachoc, P.~Moustrou, F.~Vallentin, and M.~C. Zimmermann.
\newblock Polarization of lattices: Stable cold spots and spherical designs.
\newblock \emph{arXiv preprint}, 2025.\\
\newblock \doi{10.48550/arXiv.2502.08819}.

\bibitem{Baernstein_HeatKernel_1997}
A.~{Baernstein II}.
\newblock {A minimum problem for heat kernels of flat tori}.
\newblock In \emph{{Extremal {R}iemann surfaces ({S}an {F}rancisco, {CA},
  1995)}}, volume 201 of \emph{{Contemporary Mathematics}}, page 227–243.
  American Mathematical Society, Providence, RI, 1997.\\
\newblock \doi{10.1090/conm/201/02604}.

\bibitem{RamanujanV}
B.~C. Berndt.
\newblock \emph{{Ramanujan’s Notebooks, Part V}}.
\newblock Springer, 1998.\\
\newblock \doi{10.1007/978-1-4612-1624-7}.

\bibitem{Ber28}
S.~N. Bernstein.
\newblock {Sur les fonctions absolument monotones}.
\newblock \emph{Acta Mathematica}, 52:\penalty0 1–66, 1928.
\newblock \doi{10.1007/BF02592679}.

\bibitem{BetKnu18-Born}
L.~B{\'e}termin and H.~Kn{\"u}pfer.
\newblock On {Born}'s conjecture about optimal distribution of charges for an
  infinite ionic crystal.
\newblock \emph{Journal of Nonlinear Science}, 28\penalty0 (5):\penalty0
  1629--1656, 2018.
\newblock \doi{10.1007/s00332-018-9460-3}.

\bibitem{BetFauKnu21}
L.~B\'{e}termin, M.~Faulhuber, and H.~Kn{\"u}pfer.
\newblock On the optimality of the rock-salt structure among lattices with
  charge distributions.
\newblock \emph{Mathematical Models and Methods in Applied Sciences},
  31\penalty0 (2):\penalty0 293--325, 2021{\natexlab{a}}.\\
\newblock \doi{10.1142/S021820252150007X}.

\bibitem{BetFauSte21}
L.~B\'{e}termin, M.~Faulhuber, and S.~Steinerberger.
\newblock {A variational principle for Gaussian lattice sums}.
\newblock \emph{arXiv preprint}, 2110.06008, 2021{\natexlab{b}}.
\newblock \doi{10.48550/arXiv.2110.06008}.

\bibitem{Bor22}
S.~Borodachov.
\newblock Polarization problem on a higher-dimensional sphere for a simplex.
\newblock \emph{Discrete \& Computational Geometry}, 67\penalty0 (2):\penalty0
  525--542, 2022.\\
\newblock \doi{10.1007/s00454-021-00308-1}.

\bibitem{BorHarRezSaf18}
S.~V. Borodachov, D.~P. Hardin, A.~Reznikov, and E.~B. Saff.
\newblock {Optimal discrete measures for Riesz potentials}.
\newblock \emph{{Transactions of the American Mathematical Society}},
  370:\penalty0 6973--6993, 2018.
\newblock \doi{10.1090/tran/72240}.

\bibitem{BorBor_Cubic_91}
J.~M. Borwein and P.~B. Borwein.
\newblock {A Cubic Counterpart of Jacobi's Identity and the AGM}.
\newblock \emph{Transactions of the American Mathematical Society},
  332\penalty0 (2):\penalty0 691–701, 1991.
\newblock \doi{10.2307/2001551}.

\bibitem{BoyDraHarSafSto23}
P.~G. Boyvalenkov, P.~D. Dragnev, D.~P. Hardin, E.~B. Saff, and M.~M.
  Stoyanova.
\newblock On polarization of spherical codes and designs.
\newblock \emph{{Journal of Mathematical Analysis and Applications}},
  524\penalty0 (1):\penalty0 127065, 2023.
\newblock \doi{10.1016/j.jmaa.2023.127065}.

\bibitem{CohKum07}
H.~Cohn and A.~Kumar.
\newblock {Universally optimal distribution of points on spheres}.
\newblock \emph{Journal of the American Mathematical Society}, 20\penalty0
  (1):\penalty0 99–148, 2007.\\
\newblock \doi{10.1090/S0894-0347-06-00546-7}.

\bibitem{Coh-Via22}
H.~Cohn, A.~Kumar, S.~D. Miller, D.~Radchenko, and M.~S. Viazovska.
\newblock {Universal optimality of $E_8$ and Leech lattices and interpolation
  formulas}.
\newblock \emph{Annals of Mathematics~(2)}, 196\penalty0 (3):\penalty0 983--1082,
  2022.
\newblock \doi{10.4007/annals.2022.196.3.3}.

\bibitem{ConSlo98}
J.~H. Conway and N.~J.~A. Sloane.
\newblock \emph{{Sphere Packings, Lattices and Groups}}, volume 290 of
  \emph{{Grundlehren der Mathematischen Wissenschaften}}.
\newblock Springer, New York, 3. edition, 1998.
\newblock \doi{10.1007/978-1-4757-6568-7}.

\bibitem{FauGumSha25}
M.~Faulhuber, A.~Gumber, and I.~Shafkulovska.
\newblock {The AGM of Gauss, Ramanujan's corresponding theory, and spectral
  bounds of self-adjoint operators}.
\newblock \emph{Monatshefte für Mathematik}, 206\penalty0 (3):\penalty0
  551--582, 2025.
\newblock \doi{10.1007/s00605-024-02051-0}.

\bibitem{Fek23}
M.~Fekete.
\newblock {Über die Verteilung der Wurzeln bei gewissen algebraischen
  Gleichungen mit ganzzahligen Koeffizienten}.
\newblock \emph{Mathematische Zeitschrift}, 17:\penalty0 228--249, 1923.\\
\newblock \doi{10.1007/BF01504345}.

\bibitem{Fol89}
G.~B. Folland.
\newblock \emph{{Harmonic Analysis in Phase Space}}.
\newblock Number 122 in {Annals of Mathematics Studies}. Princeton University
  Press, 1989.
\newblock \doi{10.1515/9781400882427}.

\bibitem{HarKenSaf13}
D.~P. Hardin, A.~P. Kendall, and E.~B. Saff.
\newblock {Polarization Optimality of Equally Spaced Points on the Circle for
  Discrete Potentials}.
\newblock \emph{Discrete \& Computational Geometry}, 50:\penalty0 236--243,
  2013.
\newblock \doi{10.1007/s00454-013-9502-4}.

\bibitem{HarPetSaf22}
D.~P. Hardin, M.~Petrache, and E.~B. Saff.
\newblock Unconstrained polarization ({Chebyshev}) problems: basic properties
  and {Riesz} kernel asymptotics.
\newblock \emph{Potential Analysis}, 56\penalty0 (1):\penalty0 21--64, 2022.
\newblock \doi{10.1007/s11118-020-09875-z}.

\bibitem{Mon88}
H.~L. Montgomery.
\newblock {Minimal theta functions}.
\newblock \emph{Glasgow Mathematical Journal}, 30\penalty0 (1):\penalty0
  75--85, 1988.
\newblock \doi{10.1017/S0017089500007047}.

\bibitem{PolSze31}
G.~P\'olya and G.~Szeg\H{o}.
\newblock {Über den transfiniten Durchmesser (Kapazitätskonstante) von ebenen
  und räumlichen Punktmengen}.
\newblock \emph{Journal für die reine und angewandte Mathematik},
  165:\penalty0 4--49, 1931.
\newblock \doi{10.1515/crll.1931.165.4}.

\bibitem{Wid41}
D.~V. Widder.
\newblock \emph{{The Laplace Transform}}.
\newblock Princeton University Press, 1941.\\
\newblock \doi{10.1515/9781400876457}.

\end{thebibliography}

\end{document}